\documentclass[twoside,11pt]{article}
\title{} \author{} \date{}
\usepackage{latexsym,amssymb,times}%,showkeys}
\input amssym.def
\newtheorem{te}{Theorem}[section]

\newtheorem{fac}[te]{Fact}

\newtheorem{cla}[te]{Claim}
\def\dok{\noindent{\bf Proof. }}
\def\kdok{\hfill $\Box$ \par \vspace*{2mm} }
\def\f{\varphi}
\def\o{\omega}
\def\r{\rho}
\def\s{\sigma}
\def\P{{\mathbb P}}
\def\Q{{\mathbb Q}}
\def\N{{\mathbb N}}
\def\X{{\mathbb X}}
\def\Z{{\mathbb Z}}
\def\R{{\mathbb R}}
\def\H{{\mathbb H}}
\def\K{{\mathbb K}}
\def\BG{{\mathbb G}}
\def\D{{\mathcal D}}
\def\CP{{\mathcal P}}
\def\L{{\mathcal L}}
\def\la{\langle}
\def\ra{\rangle}

\def\supp{\mathop{\mbox{supp}}\nolimits}
\def\Rado{\mathop{\mathrm{Rado}}\nolimits}
\begin{document}
\thispagestyle{plain}
\begin{center}
           {\large \bf
           MAXIMAL CHAINS OF ISOMORPHIC SUBGRAPHS\\[1mm]
           OF COUNTABLE ULTRAHOMOGENEOUS GRAPHS}
\end{center}
\begin{center}
{\bf Milo\v s S.\ Kurili\'c\footnote{Department of Mathematics and
Informatics, University of Novi Sad,
     Trg Dositeja Obradovi\'ca 4, 21000 Novi Sad, Serbia. e-mail: milos@dmi.uns.ac.rs}
     and
     Bori\v sa  Kuzeljevi\'c\footnote{Mathematical Institute of the Serbian Academy of Sciences and Arts,
     Kneza Mihaila 36, 11001 Belgrade, Serbia. e-mail: borisa@mi.sanu.ac.rs}
}
\end{center}
\begin{abstract}
\noindent
For a countable ultrahomogeneous graph $\BG=\la G, \r\ra$ let $\P (\BG)$ denote the collection of sets $A\subset G$ such that
$\la A, \r \cap [A]^2\ra \cong \BG$.  The order types of maximal chains in the poset
$\langle \P(\BG )\cup \{\emptyset \}, \subset \rangle$ are characterized as:

(I) the order types of compact sets of reals having
the minimum non-isolated, if $\BG$ is the Rado graph or the Henson graph $\H _n$, for some $n\geq 3$;

(II) the order types of compact nowhere dense sets of reals having
the minimum non-isolated, if $\BG$ is the union of $\mu$ disjoint complete graphs of size $\nu$, where $\mu \nu =\o$.

\end{abstract}
\section{Introduction}\label{S0}
If $\X$ is a relational structure, $\P (\X )$ will denote the set of domains of substructures of $\X$ which are isomorphic to $\X$.
$\X$ is called {\it ultrahomogeneous} iff each isomorphism between two finite substructures of $\X$ can be extended to an automorphism of $\X$.

A structure $\BG =\la G,\r \ra$ is a {\it graph} iff $G$ is a set and $\r$ a symmetric irreflexive binary relation on $G$.
We will also use the following equivalent definition: a pair $\BG =\la G,\r \ra$ is a graph iff $G$ is a set and $\r \subset [G]^2$.
Then for $H\subset G$, $\la H , \r \cap [H]^2\ra $
(or $\la H , \r \cap (H\times H)\ra $, in the relational version) is the corresponding {\it subgraph} of $\BG$.
For a cardinal $\nu$,  $\K _\nu$ will denote the {\it complete graph} of size $\nu$.
A graph is called $\K _n$-free iff it has no subgraphs isomorphic to $\K _n$.
We will use the following well-known classification of
countable ultrahomogeneous graphs \cite{Lach}:
\begin{te}[Lachlan and Woodrow]\rm\label{T6039}
Each countable ultrahomogeneous graph is isomorphic to one of the following graphs

- $\BG _{\mu \nu }$, the union of $\mu$ disjoint copies of $\K _\nu$, where $\mu \nu=\o$,

- $\BG_{\Rado}$, the unique countable homogeneous universal graph, the Rado graph,

- $\H _n$, the unique countable homogeneous universal $\K _n$-free graph, for $n\geq 3$,

- the complements of these graphs.
\end{te}
Properties of maximal chains in posets are widely studied order invariants (see \cite{Day}, \cite{Kopp}, \cite{Kura}, \cite{McMo}, \cite{Monk})
and, as a part of investigation of the partial orders of the form $\la \P (\X ),\subset \ra$, where $\X$ is a relational structure,
the class of order types of maximal chains in the poset $\la \P (\BG_{\Rado}),\subset \ra$
was characterized in \cite{KK}. The aim of this paper is to complete the picture for all countable ultrahomogeneous graphs in this context and, thus,
the following theorem is our main result.
\begin{te}\rm\label{T6021}
Let $\BG$ be a countable ultrahomogeneous graph. Then

\noindent
(I) If $\BG = \BG_{\Rado}$ or $\BG =\H _n$, for some $n\geq3$, then for each linear order $L$ the
following conditions are equivalent:

(a) $L$ is isomorphic to a maximal chain in the poset $\langle \P(\BG ) \cup \{ \emptyset \} , \subset \rangle $;

(b) $L$ is an ${\mathbb R}$-embeddable complete linear order with $0_L$ non-isolated;

(c) $L$ is isomorphic to a compact set $K \subset \mathbb R$ having the minimum non-isolated.

\noindent
(II) If $\BG = \BG_{\mu \nu}$, where $\mu \nu =\o$, then for each linear order $L$ the
following conditions are equivalent:

(a) $L$ is isomorphic to a maximal chain in the poset $\langle \P(\BG ) \cup \{ \emptyset \} , \subset \rangle $;

(b) $L$ is an ${\mathbb R}$-embeddable Boolean linear order with $0_L$ non-isolated;

(c) $L$ is isomorphic to a compact nowhere dense set $K \subset \mathbb R$ having the minimum non-isolated.
\end{te}
It is easy to check that for a relational structure $\la X, \r\ra$ we have $\P (\la X, \r\ra)=\P (\la X, \r ^c\ra)$ and, hence,
regarding Theorem \ref{T6039}, Theorem \ref{T6021} in fact covers all countable ultrahomogeneous graphs.

The statement (I) for the Rado graph is proved in \cite{KK} and a proof for the graphs $\H _n$ is given in Section \ref{S4},
while (II) is proved in Section \ref{S5}.
\section{Preliminaries}\label{S1}
In this section we recall basic definitions and facts which will be used in the paper.

If $\langle P,\leq
\rangle$ is a partial order, then the {\it smallest} and the {\it
largest element} of $P$ are denoted by $0_P$ and $1_P$; the {\it
intervals} $(x,y)_P$, $[x,y]_P$, $(-\infty , x)_P$ etc.\ are
defined in the usual way. A set $D\subset P$ is {\it dense} iff
for each $p\in P$ there is $q\in D$ such that $q\leq p$. A set $G\subset P$ is a {\it filter} iff (F1) for each $p,q \in G$ there is $r\in G$
such that $r\leq p, q$ and (F2) $G\ni p\leq q$ implies $q\in G$.
\begin{fac}[Rasiowa-Sikorski] \rm \label{T6023}
If $D_n$, $n\in \omega$ are dense sets in a partial order $\langle
P, \leq \rangle$, then there is a filter $G$ in $P$ intersecting
all of them.
\end{fac}
\dok
Let $p_0 \in D_0$ and, for $n\in \omega$, let us pick $
p_{n+1}\in D_{n+1} $ such that $ p_{n+1} \leq p_n$. Then $G=\{
p\in P: \exists n\in \omega \; p_n \leq p\}$ is a filter
intersecting all $D_n$'s.
\kdok
A pair $\langle{\mathcal A},{\mathcal B}\rangle$ is a {\it cut} in
a linear order $\langle L,<\rangle$ iff $L={\mathcal
A}\;{\buildrel {.} \over {\cup}}\; {\mathcal B}$, ${\mathcal
A},{\mathcal B}\neq \emptyset$ and $a<b$, for each $a\in {\mathcal
A}$ and $b\in {\mathcal  B}$. A cut $\langle {\mathcal
A},{\mathcal B} \rangle$ is  a {\it gap} iff neither $\max
{\mathcal A}$ nor $\min {\mathcal  B}$ exist.
$\langle L,<\rangle $ is called:
{\it complete} iff  it has 0 and 1 and has no gaps; {\it dense}
iff $(x,y)_L\neq \emptyset$, for each $x,y \in L$ satisfying
$x<y$; {\it ${\mathbb R}$-embeddable} iff it is isomorphic to a
subset of ${\mathbb R}$; {\it Boolean} iff it is complete and {\it
has dense jumps}, which means that for each $x,y \in L$ satisfying
$x<y$ there are $a,b \in L$ such that $x\leq a< b \leq y$ and
$(a,b)_L=\emptyset$. A set $D \subset L$  is called {\it dense in}
$L$ iff for each $x,y \in L$ satisfying $x<y$ there is $z\in D$
such that $x< z < y$. If $\langle I , <_I \rangle $ and  $\langle
L_i,<_i\rangle$, $i\in I$, are linear orders and $L_i\cap
L_j=\emptyset$,  whenever $i\neq j$,  then the corresponding {\it
lexicographic sum} $\sum_{i\in I}L_i$ is the linear ordering
$\langle \bigcup_{i\in I}L_i,<\rangle$ where the relation $<$ is defined by: $x<y  \Leftrightarrow \exists
i\in I \; (x,y\in L_i \wedge x<_iy)
 \vee \exists i,j\in I \; (i<_I j\wedge x\in L_i \wedge y\in L_j)$.
\begin{fac} \rm \label{T6024}
If $\langle L, <\rangle$ is an at most countable complete linear
order, it is Boolean.
\end{fac}
\dok
Let $x,y\in L$ and $x<y$. Suppose that for each $a,b \in [x,y
]_L$ satisfying $a<b$ we have $(a,b)_L \neq \emptyset$. Then $[x,y
]_L$ would be a dense complete linear order, which is impossible
because $L$ is countable. Thus $L$ has dense jumps.
\kdok
\noindent
${\mathcal P} \subset P(\omega)$ is called a {\it positive family} iff
(P1) $\emptyset \notin {\mathcal P}$;
(P2) ${\mathcal P} \ni A \subset B \subset \omega \Rightarrow  B
\in {\mathcal P}$;
(P3) $A \in {\mathcal P} \wedge |F| <\omega \Rightarrow  A
\backslash F \in {\mathcal P}$;
(P4) $\exists A \in{\mathcal P} \; |\omega \backslash
A|=\omega$.
\begin{fac} \rm \label{T6022}
(See \cite{K1})
If ${\mathcal P}\subset P(\omega)$ is a positive family, then for each linear order $L$ the
following conditions are equivalent:

(a) $L$ is isomorphic to a maximal chain in the poset $\langle \CP \cup \{ \emptyset \} , \subset \rangle $;

(b) $L$ is an ${\mathbb R}$-embeddable Boolean linear order with $0_L$ non-isolated;

(c) $L$ is isomorphic to a compact nowhere dense set $K \subset \mathbb R$ having the minimum non-isolated.

(d) $L$ is isomorphic to a maximal chain $\mathcal L$ in the poset $\langle {\mathcal P} \cup \{\emptyset\}, \subset \rangle$
such that $\bigcap (\mathcal L \setminus \{ \emptyset \}) =\emptyset$.
\end{fac}
\begin{fac}  \rm \label{T6030}
Let $A\subset B\subset \omega$ and let $L$ be a complete linear
ordering, such that $|B\setminus A|=|L|-1$. Then there is a chain
${\mathcal L}$ in $[A,B]_{P (B )}$ satisfying $A,B\in {\mathcal
L}\cong L$ and such that $\bigcup {\mathcal  A},\bigcap {\mathcal
B}\in {\mathcal  L}$ and $|\bigcap {\mathcal  B}\setminus \bigcup
{\mathcal  A}|\leq 1$, for each cut $\langle {\mathcal
A},{\mathcal B}\rangle$  in ${\mathcal L}$.
\end{fac}
\dok
If $|B\setminus A|$ is a finite set, say $B=A \cup \{ a_1 ,
\dots a_n \}$, then $|L|=n+1$ and $\mathcal L =\{ A, A\cup \{ a_1
\}, A\cup \{ a_1 , a_2  \} , \dots , B\}$ is a chain with the
desired properties.

If $|B\setminus A|=\omega$, then $L$ is a countable and, hence,
$\mathbb R$-embeddable complete linear order. It is known that an
infinite linear order is isomorphic to a maximal chain in
$P(\omega)$ iff it is ${\mathbb R}$-embeddable and Boolean (see,
for example, \cite{K}). By Fact \ref{T6024} $L$ is a Boolean order
and, thus, there is a maximal chain $\mathcal L _1$ in
$P(B\setminus A)$ isomorphic to $L$. Let $\mathcal L =\{ A\cup C :
C\in \mathcal L _1\}$. Since $\emptyset , B\setminus A \in
\mathcal L_1$ we have $A,B \in \mathcal L$ and the function $f:
\mathcal L _1 \rightarrow \mathcal L$, defined by $f(C)=A\cup C$,
witnesses that $\langle \mathcal L _1 , \varsubsetneq \rangle
\cong \langle \mathcal L , \varsubsetneq \rangle$ so $\mathcal L$
is isomorphic to $L$. For each cut $\langle \mathcal A , \mathcal
B\rangle$ in $\mathcal L_1$ we have $\bigcup \mathcal A \subset
\bigcap \mathcal B$ and, by the maximality of $\mathcal L_1$,
$\bigcup \mathcal A , \bigcap \mathcal B \in \mathcal L _1$  and
$| \bigcap \mathcal B \setminus \bigcup \mathcal A | \leq 1$.
Clearly, the same is true for each cut in $\mathcal L$.
\hfill $\Box$
\section{General results}\label{S2}
The following three general statements, concerning the class of the order types of
maximal chains of copies of relational structures, will be used in
the proof of Theorem \ref{T6021}. The first one gives a necessary condition for a type to be in the class
corresponding to a countable ultrahomogeneous structure.
\begin{te} \rm \label{T5037}
(\cite{KK1}) Let $\mathbb X$ be a countable ultrahomogeneous structure of an at
most countable relational language and $\mathbb P (\mathbb X )\neq
\{ X\}$. Then for each linear order $L$ we have (a) $\Rightarrow$ (b), where

(a) $L$ is isomorphic to a maximal chain in the poset $\langle \P(\X ) \cup \{ \emptyset \} , \subset \rangle $;

(b) $L$ is an ${\mathbb R}$-embeddable complete linear order with $0_L$ non-isolated.

\noindent
In particular, the union of a chain in $\langle \P(\X ) , \subset \rangle $ belongs to $\langle \P(\X )$.
\end{te}
The following statement describes a class of structures such that, regarding Theorem \ref{T5037},
the implication (b) $\Rightarrow$ (a) holds for each linear order $L$.
\begin{te}\rm \label{T6037}
Let $\X$ be a countable relational structure and $\Q$ the set of rationals.
\begin{itemize}
\item[(A)] If there exist a partition $\{ J_n:n\in \o \}$ of $\Q$ and a structure with the domain $\Q$ of the same signature as $\X$ such that
\begin{itemize}
\item[(i)] $J_0$ is a dense subset of $\Q$,

\item[(ii)] $J_n$, $n\in \N$, are coinitial subsets of $\Q$,

\item[(iii)] $J_0 \cap(-\infty ,x) \subset  A\subset \Q\cap (-\infty , x)$ implies $A\cong \X$, for all $x \in \R \cup \{ \infty \}, $

\item[(iv)]  $J_0 \cap(-\infty ,q] \subset  C\subset \Q \cap (-\infty , q] $ implies $C\not\cong \X$, for each $q\in J_0$,
\end{itemize}
then for each uncountable $\mathbb R$-embeddable complete linear order $L$ with $0_L$ non-isolated and such that
           all initial segments of $L\setminus \{ 0_L\}$ are uncountable there is a maximal chain in
           $\langle \mathbb P(\X)\cup \{\emptyset\},\subset \rangle$ isomorphic to $L$.

\item[(B)] If, in addition,
\begin{itemize}
\item[(v)] for each countable complete linear order $L$ with $0_L$ non-isolated
there is a maximal chain in $\langle \mathbb P(\X)\cup \{\emptyset\},\subset \rangle$ isomorphic to $L$,
\end{itemize}
\noindent
then for each $\mathbb R$-embeddable complete linear order $L$ with $0_L$ non-isolated
        there is a maximal chain in $\langle \mathbb P(\X)\cup \{\emptyset\},\subset \rangle$ isomorphic to $L$.
\end{itemize}
\end{te}
\dok
Let $L$ be an uncountable $\mathbb R$-embeddable complete linear order with $0_L$ non-isolated.

\begin{cla}  \rm \label{T6014}
$L\cong \sum_{x \in [-\infty , \infty ]} L_x$, where

(L1) $L_x$, $x\in [-\infty , \infty ]$, are at most countable
complete linear orders,

(L2) The set $M=\{x \in [-\infty , \infty ]:|L_x| >1 \}$ is at
most countable,

(L3) $|L_{-\infty }|=1$ or $0_{L_{-\infty }}$ is non-isolated.
\end{cla}
\dok
$L=\sum _{i\in I}L_i$, where $L_i$ are the equivalence
classes corresponding to the condensation relation $\sim$ on $L$
given by: $x\sim y \Leftrightarrow |[\min \{ x,y \} ,\max \{ x,y
\} ]|\leq \omega$ (see \cite{Rosen}). Since $L$ is complete and
${\mathbb R}$-embeddable $I$ is too and, since the cofinalities
and coinitialities of $L_i$'s are countable, $I$ is a dense linear
order; so $I\cong [0,1]\cong [-\infty , \infty ]$. Hence $L_i$'s
are complete and, since $\min L_i \sim \max L_i$, countable. If
$|L_i|>1$,  $L_i$ has a jump (Fact \ref{T6024}) so,
$L\hookrightarrow {\mathbb R}$ gives $|M|\leq \omega$.
\kdok

(A) Let all initial segments of $L\setminus \{ 0_L\}$ be uncountable. Then, by Claim \ref{T6014},  $|L_{-\infty }|=1$,  that is $-\infty\notin M$, and we have two cases.

\noindent
{\bf Case I:} $\infty \in M$.
By (L2) there is an injection $\f : M\rightarrow \N$.
By (L1), for $y \in M$ we have $|L_y|\leq \o $ and by (ii) $|J_{\f (y)} \cap (-\infty,y)|=\o$ so we take $I_y \in [J_{\f (y)} \cap (-\infty,y)]^{|L_y|-1}$.
Let us define the sets $A_x$, $x\in [-\infty , \infty]$ and $A_x^+$, $x\in M$, by
$$ \textstyle
A_x= \left\{\begin{array}{ll}
\emptyset                  ,                                                   & \mbox{for }  x = -\infty, \\
\big(J_0 \cap (-\infty,x)  \big)\cup \bigcup{_{y \in M \cap (-\infty,x)}}I_y , & \mbox{for }  x \in (-\infty,\infty];
\end{array}
\right.
$$
\vspace{-2mm}
$$\textstyle
A_x^+=A_x\cup I_x, \quad \mbox{for } x \in M.
$$
Since $J_0 \subset A_\infty ^+  = J_0 \cup \bigcup _{y\in M}I_y \subset \mathbb Q$, by (iii) we have $A_\infty ^+ \cong \X$ and we construct a maximal chain
${\mathcal L}$ in $\langle \mathbb P(A_\infty ^+ )\cup \{\emptyset \}, \subset\rangle$, such that ${\mathcal L}\cong L$.
\begin{cla}  \rm \label{T6031}
The sets  $ A_x $, $x\in [-\infty,\infty]$ and $A_x^{+}$, $x \in
M$ are subsets of the set $A_{\infty}^+$. In
addition, for each $x, x_1, x_2 \in [-\infty,\infty]$ we have

(a) $A_x \subset (-\infty,x)$;

(b) $A_x^+ \subset (-\infty,x)$, if $x \in M$;

(c) $x_1 < x_2 \Rightarrow A_{x_1}\varsubsetneq A_{x_2}$;

(d) $M \ni x_1 < x_2 \Rightarrow  A_{x_1}^+ \varsubsetneq
A_{x_2}$;

(e)  $|A_x^+\setminus A_x|=|L_x|-1$, if $x \in M$;

(f) $ A_x \in \mathbb P (A_{\infty}^+ )$, for each $x\in (-\infty
, \infty ]$.

(g) $ A_x ^+ \in \mathbb P (A_{\infty}^+  )$ and $[A_x , A_x ^+]_{\mathbb P (A_{\infty}^+  )} = [A_x , A_x ^+]_{P( A_x ^+ )}$,
for each $x\in M$.
\end{cla}
\dok
Statements (c) and (d) are true since $J_0$ is a dense subset
of $\mathbb Q$; (a), (b) and (e) follow from the definitions of
$A_x$ and $A_x^+$ and the choice of the sets $I_y$.
For  $x\in (-\infty , \infty ]$ we have $J_0 \cap (-\infty,x) \subset A_x \subset \Q \cap (-\infty,x)$ so, by (iii),
$A_x \cong \X \cong A_{\infty}^+$ and (f) is true.
If $x\in M$, then $J_0 \cap (-\infty,x) \subset A_x \subset A_x^+ \subset\Q \cap (-\infty,x)$ so, by (iii),
$A_x \subset A \subset A_x^+$ implies $A\cong \X\cong A_{\infty}^+$
 and (g) is true as well.
\kdok

\noindent Now, for $x \in [-\infty,\infty]$ we define chains
${\mathcal L}_x $ in $\la \mathbb P (A_{\infty}^+  )\cup \{\emptyset \}, \subset \ra$ as follows.

For $x \not\in M$ we define ${\mathcal L}_x=\{A_x\}$. In
particular, ${\mathcal L}_{-\infty }=\{ \emptyset \}$.

For $x \in M$, using Claim \ref{T6031}(g) and Fact \ref{T6030} we
obtain ${\mathcal L}_x \subset [A_x , A_x ^+]_{P(A_x ^+ )}$
such that $\langle \mathcal L _x , \varsubsetneq \rangle \cong
\langle L_x , <_x \rangle $ and
\begin{equation}\label{EQ6012}\textstyle
A_x, A_x^+ \in {\mathcal L}_x \subset [A_x,A_x^+]_{\mathbb P
(A_{\infty}^+ )},
\end{equation}
\begin{equation}\label{EQ6013}\textstyle
\bigcup {\mathcal  A},\bigcap {\mathcal  B}\in {\mathcal  L}_x
\;\;\mbox{ and }\;\; |\bigcap {\mathcal  B}\setminus \bigcup
{\mathcal  A}|\leq 1, \mbox{ for each cut }\langle {\mathcal
A},{\mathcal B}\rangle \mbox{ in }{\mathcal L}_x.
\end{equation}
For ${\mathcal A},{\mathcal B} \subset \mathbb P (A_{\infty}^+)$
we will write ${\mathcal A}\prec {\mathcal B}$ iff $A\varsubsetneq
B$, for each $A \in {\mathcal A}$ and $B \in {\mathcal B}$.
\begin{cla}  \rm \label{T6033}
Let ${\mathcal L}=\bigcup_{x \in [-\infty,\infty]}{\mathcal L}_x$.
Then

(a) If $-\infty \leq x_1 <x_2  \leq \infty $, then ${\mathcal
L}_{x_1}\prec {\mathcal L}_{x_2}$ and $\bigcup{\mathcal L}_{x_1}
\subset A_{x_2}\subset \bigcup{\mathcal L}_{x_2}.$

(b) ${\mathcal L}$ is a chain in $\langle \mathbb P
(A_{\infty}^+)\cup \{ \emptyset \},\subset\rangle$ isomorphic to
$L= \sum_{x \in [-\infty,\infty]}L_x$.

(c) ${\mathcal L}$ is a maximal chain in $\langle \mathbb P
(A_{\infty}^+ )\cup \{ \emptyset \},\subset\rangle$.
\end{cla}
\dok
(a) Let $A \in {\mathcal L}_{x_1}$ and $B \in {\mathcal L}_{x_2}$. If $x_1 \in (-\infty,\infty]\setminus M$, then, by
(\ref{EQ6012}) and Claim \ref{T6031}(c) we have
$A=A_{x_1}\varsubsetneq A_{x_2} \subset B$. If $x_1 \in M$, then,
by (\ref{EQ6012}) and Claim \ref{T6031}(d), $A \subset
A_{x_1}^+\varsubsetneq A_{x_2} \subset B$. The second statement
follows from $A_{x_2} \in {\mathcal L}_{x_2}$.

(b) By (a), $\langle [-\infty,\infty],<\rangle \cong \langle
\{{\mathcal L}_x:x \in [-\infty,\infty]\},\prec \rangle$. Since
${\mathcal L}_x \cong L_x$, for $x \in [-\infty,\infty]$, we have
$\langle {\mathcal L},\varsubsetneq\rangle \cong\sum_{x \in
[-\infty,\infty]}\langle {\mathcal L}_x,\varsubsetneq\rangle \cong
\sum_{x \in [-\infty,\infty]}L_x =L$.

(c) Suppose that $C\in \mathbb P (A_{\infty}^+ )\cup \{ \emptyset
\}$ witnesses that ${\mathcal L}$ is not maximal. Clearly
${\mathcal L}={\mathcal A}\dot{\cup } {\mathcal B}$ and ${\mathcal
A}\prec {\mathcal B}$, where ${\mathcal A}=\{ A \in {\mathcal L}:
A \varsubsetneq C \}$ and  ${\mathcal B}=\{ B \in {\mathcal L}: C
\varsubsetneq B \}$. Now $\emptyset \in {\mathcal L}_{-\infty }$
and, since $\infty \in M$, by (\ref{EQ6012}) we have $A_{\infty}^+
\in {\mathcal L}_{\infty}$. Thus $\emptyset ,A_{\infty}^+ \in
{\mathcal L}$, which implies ${\mathcal A},{\mathcal B}\neq
\emptyset$ and, hence, $\langle {\mathcal A},{\mathcal B}\rangle$
is a cut in $\langle {\mathcal L},\varsubsetneq\rangle$. By
(\ref{EQ6012}) we have  $\{ A_x :x \in (-\infty,\infty ]\}\subset
{\mathcal L}\setminus \{ \emptyset \}$ and, by Claim
\ref{T6031}(a), $\bigcap ({\mathcal L}\setminus \{ \emptyset
\})\subset \bigcap_{x \in (-\infty,\infty]}A_x \subset \bigcap
_{x\in (-\infty,\infty]}(-\infty,x)=\emptyset$, which implies
${\mathcal A}\neq \{ \emptyset\}$. Clearly,
\begin{equation}\label{EQ6014}\textstyle
\bigcup{\mathcal A} \subset C \subset \bigcap {\mathcal B}.
\end{equation}
\noindent {\it Case 1:} ${\mathcal A}\cap {\mathcal
L}_{x_0}\neq\emptyset$ and ${\mathcal B}\cap {\mathcal
L}_{x_0}\neq \emptyset$, for some $x_0 \in (-\infty,\infty]$. Then
$|{\mathcal L}_{x_0}|>1$, $x_0 \in M$ and $\langle {\mathcal
A}\cap{\mathcal L}_{x_0},{\mathcal B}\cap {\mathcal
L}_{x_0}\rangle$ is a cut in ${\mathcal L}_{x_0}$ satisfying
(\ref{EQ6013}). By (a), ${\mathcal A}=\bigcup_{x<x_0}{\mathcal
L}_x\cup ({\mathcal A}\cap {\mathcal L}_{x_0})$ and, consequently,
$\bigcup{\mathcal A}=\bigcup({\mathcal A}\cap{\mathcal L}_{x_0})
\in {\mathcal L}$. Similarly, $\bigcap {\mathcal
B}=\bigcap({\mathcal B}\cap{\mathcal L}_{x_0})\in {\mathcal L}$
and, since $|\bigcap {\mathcal B}\setminus \bigcup{\mathcal
A}|\leq 1$, by (\ref{EQ6014}) we have $C\in {\mathcal L}$. A
contradiction.

\vspace{2mm} \noindent {\it Case 2:} $\neg$ Case 1. Then for each
$x \in (-\infty,\infty]$ we have ${\mathcal L}_x \subset {\mathcal
A}$ or ${\mathcal L}_x \subset {\mathcal B}$. Since ${\mathcal
L}={\mathcal A}\;{\buildrel {.}\over{\cup}}\; {\mathcal B}$,
${\mathcal A}\neq \{ \emptyset\}$ and ${\mathcal A}, {\mathcal
B}\neq\emptyset$, the sets ${\mathcal A}'=\{x \in
(-\infty,\infty]:{\mathcal L}_x\subset {\mathcal A}\}$ and
${\mathcal B}'=\{x \in(-\infty,\infty]:{\mathcal L}_x \subset
{\mathcal B}\}$ are non-empty and $(-\infty,\infty]={\mathcal
A}'\;{\buildrel {.} \over {\cup}}\;{\mathcal B}'$. Since
${\mathcal A}\prec {\mathcal B}$, for $x_1 \in {\mathcal A}'$ and
$x_2 \in {\mathcal B}'$ we have ${\mathcal L}_{x_1}\prec {\mathcal
L}_{x_2}$ so, by (a), $x_1<x_2$. Thus $\langle {\mathcal
A}',{\mathcal B}' \rangle$ is a cut in $(-\infty,\infty]$ and,
consequently, there is $x_0 \in (-\infty,\infty]$ such that
$x_0=\max {\mathcal A}'$ or $x_0=\min {\mathcal B}'$.

\vspace{2mm} \noindent {\it Subcase 2.1:} $x_0=\max {\mathcal
A}'$. Then $x_0<\infty$ because $\mathcal B \neq \emptyset$ and
${\mathcal A}=\bigcup_{x\leq x_0}{\mathcal L}_x$ so, by (a),
$\bigcup{\mathcal A}=\bigcup_{x\leq x_0}\bigcup{\mathcal
L}_x=\bigcup_{x<x_0}\bigcup{\mathcal L}_x \cup \bigcup{\mathcal
L}_{x_0}=\bigcup{\mathcal L}_{x_0}$ which, together with
(\ref{EQ6012}) implies
\begin{equation}\label{EQ6015}\textstyle
\bigcup{\mathcal A}=\left\{\begin{array}{ll}
                         A_{x_0} & \mbox{if } x_0 \not\in M, \\
                         A_{x_0}^+ & \mbox{if } x_0 \in M.
                         \end{array}
                         \right.
\end{equation}
%$$
Since ${\mathcal B}=\bigcup_{x\in (x_0,\infty]}{\mathcal L}_x$, we
have
$\bigcap {\mathcal B}=\bigcap_{x \in(x_0,\infty]}\bigcap{\mathcal L}_x$. By (\ref{EQ6012})
$\bigcap{\mathcal L}_x=A_x$,
so we have
$\bigcap {\mathcal B}
=\big(\bigcap{_{x \in (x_0,\infty]}}(-\infty,x)\cap J_0\big) \cup \big(\bigcap{_{x \in (x_0,\infty]}}\bigcup{_{y \in M\cap(-\infty,x)}}I_y\big)
=\big((-\infty,x_0]\cap J_0 \big) \cup \bigcup{_{y \in M \cap(-\infty,x_0]}}I_y
=\textstyle A_{x_0} \cup \big(\{x_0\}\cap J_0 \big) \cup \bigcup{_{y \in M\cap\{x_0\}}}I_y ,$
so
\begin{equation}\label{EQ6016}\textstyle
\bigcap{\mathcal B}=\left\{\begin{array}{llr}
                          A_{x_0} & \mbox{if} & x_0 \notin J_0 \quad \wedge \quad x_0 \notin M, \\
                          A_{x_0}\cup\{x_0\} & \mbox{if} & x_0 \in J_0 \quad \wedge \quad x_0 \notin M, \\
                          A_{x_0}^+ & \mbox{if} & x_0 \notin J_0 \quad \wedge \quad x_0 \in M, \\
                          A_{x_0}^+ \cup \{x_0\} & \mbox{if} & x_0 \in J_0 \quad \wedge \quad x_0 \in M.
                          \end{array}
                          \right.
\end{equation}
%$$
If $x_0 \not\in J_0$, then, by (\ref{EQ6014}), (\ref{EQ6015}) and
(\ref{EQ6016}), we have $ \bigcup \mathcal A =\bigcap \mathcal B
=C \in \mathcal L$. A contradiction.

If $x_0 \in J_0$ and $x_0 \not\in M$, then $\bigcup \mathcal A = A_{x_0}$ and
$\bigcap \mathcal B = A_{x_0} \cup \{ x_0\}$. So, by
(\ref{EQ6014}) and since $C\not\in \mathcal L$ we have
$C=\bigcap \mathcal B= A_{x_0} \cup \{ x_0\}$.
Thus $J_0 \cap (-\infty , x_0]\subset C$ and, by Claim \ref{T6031}(a), $C\subset (-\infty , x_0]$.
But by (iv) we have $C\not\cong \X (\cong A_\infty ^+)$.
A contradiction.

If $x_0 \in J$ and $x_0 \in M$, then $\bigcup \mathcal A =
A_{x_0}^+$ and $\bigcap \mathcal B = A_{x_0}^+ \cup \{ x_0\}$.
Again, by (\ref{EQ6014}) and since $C\not\in \mathcal L$ we have
$C=\bigcap \mathcal B= A_{x_0}^+ \cup \{ x_0\}$.
Thus $J_0 \cap (-\infty , x_0]\subset C$ and, by Claim \ref{T6031}(b), $C\subset (-\infty , x_0]$.
Again, by (iv) we have $C\not\cong \X (\cong A_\infty ^+)$, a contradiction.

\vspace{2mm} \noindent {\it Subcase 2.2:} $x_0=\min {\mathcal
B}'$. Then, by (\ref{EQ6012}), $A_{x_0} \in {\mathcal L}_{x_0}
\subset {\mathcal B}$ which, by (a), implies $\bigcap{\mathcal
B}=A_{x_0}$. Since $A_x \in {\mathcal L}_x$, for $x \in
(-\infty,\infty]$ and ${\mathcal A}=\bigcup_{x<x_0}{\mathcal L}_x$
we have $\bigcup{\mathcal A} =\bigcup_{x<x_0}\bigcup{\mathcal L}_x
\supset \bigcup_{x<x_0}A_x =\bigcup_{x<x_0}\big((-\infty,x)\cap J_0
\big)\cup \bigcup_{x<x_0}\bigcup_{y \in M \cap(-\infty,x)}I_y$
$=\big((-\infty,x_0)\cap J_0 \big)\cup \bigcup_{y \in M
\cap(-\infty,x_0)}I_y =A_{x_0}$ so $A_{x_0} \subset
\bigcup{\mathcal A}\subset \bigcap{\mathcal B} = A_{x_0}$, which
implies $C=A_{x_0} \in {\mathcal L}$. A contradiction.
\kdok

\noindent {\bf Case II:} $\infty \not\in M $. Then
$L_\infty =\{ \max L\}$ and the sum $L+1$ belongs to Case I. So,
there exists a maximal chain ${\mathcal L}$ in $\langle \mathbb P(\X )\cup \{\emptyset \}, \subset \rangle $
and an isomorphism $f: \langle L+1 , < \rangle \rightarrow \langle {\mathcal L} , \subset\rangle $.
Then $A=f(\max L)\in \mathbb P (\X )$ and ${\mathcal L}'=f[L]\cong L$. By the maximality of ${\mathcal L}$,
${\mathcal L}'$ is a maximal chain in $\langle \mathbb P (A)\cup \{\emptyset \}, \subset \rangle \cong \langle \mathbb P (\X )\cup \{\emptyset \}, \subset \rangle$.

\vspace{2mm}

(B) Since (v) holds we assume that $L$ is uncountable.
If all initial segments of $L$ are uncountable, the statement is proved in (A). Otherwise, by Claim \ref{T6014}
we have $L=\sum_{x \in [-\infty,\infty]}L_x$, (L1) and (L2) hold and
\begin{center}
(L3$'$) $L_{-\infty}$ is a countable complete linear order with
$0_{L_{-\infty}}$ non-isolated.
\end{center}
\noindent Clearly $L=L_{-\infty }+ L^+$, where $L^+ =\sum_{x \in
(-\infty,\infty]}L_x = \sum_{y \in (0,\infty]}L_{\ln y}$ (here
$\ln \infty =\infty$). Let $L_y'$, $y \in [-\infty,\infty]$, be
disjoint linear orders such that $L_y'\cong 1$, for $y \in
[-\infty , 0]$, and $L_y'\cong L_{\ln y}$, for $y\in (0, \infty
]$. Now $\sum_{y \in [-\infty,\infty]}L_y' \cong [-\infty , 0] +
L^+$ and by (A) we obtain a maximal chain
${\mathcal L}$ in $\mathbb P (\X )\cup \{ \emptyset \}$ and an
isomorphism $f: \langle [-\infty , 0] + L^+ , < \rangle
\rightarrow \langle {\mathcal L}, \subset \rangle$. Clearly, for
$A_0=f(0)$ and ${\mathcal L}^+=f[L^+]$ we have $A_0\in {\mathcal
L}$ and ${\mathcal L}^+\cong L^+$.

By the assumption and (L3$'$), $\mathbb P (A_0)\cup \{
\emptyset \}$ contains a maximal chain ${\mathcal
L}_{-\infty}\cong L_{-\infty }$. Clearly $A_0\in {\mathcal
L}_{-\infty}$ and ${\mathcal L}_{-\infty} \cup {\mathcal L}^+
\cong L_{-\infty }+L^+ =L$. Suppose that $B$ witnesses that
${\mathcal L}_{-\infty} \cup {\mathcal L}^+ $ is not a maximal
chain in $\mathbb P (\X )\cup \{ \emptyset \}$. Then either $A_0
\varsubsetneq B$, which is impossible since ${\mathcal L}$ is
maximal in $\mathbb P (\X )\cup \{ \emptyset \}$, or
$B\varsubsetneq A_0$, which is impossible since ${\mathcal
L}_{-\infty}$ is maximal in $\mathbb P(A_0)\cup \{ \emptyset \}$.
\kdok
The following theorem gives a sufficient condition for (v) of Theorem \ref{T6037}.
\begin{te}\rm\label{T6025}
Let $\X =\la X, \la \s _i :i\in I \ra\ra$ be a countable relational structure.
If there is a positive family $\CP$ in $P(X)$ such that $\CP \subset \P (\X )$ and $\bigcap \CP =\emptyset$,
then

(a) For each ${\mathbb R}$-embeddable Boolean linear order $L$ with $0_L$ non-isolated there is a maximal chain in
$\la \P(\X )\cup\{\emptyset\}, \subset\ra$ isomorphic to $L$;

(b) For each countable complete linear order $L$ with $0_L$ non-isolated there is a maximal chain in
$\la \P(\X )\cup\{\emptyset\}, \subset\ra$ isomorphic to $L$.
\end{te}
\dok
(a) By Fact \ref{T6022} there is a maximal chain $\L$ in $\CP\cup\{\emptyset\}$ isomorphic
to $L$ and satisfying $\bigcap (\mathcal L\setminus\{\emptyset\})=\emptyset$.
Suppose that $C\in \P(\X )\cup \{\emptyset \}$ witnesses that $\L$ is not a maximal chain in
$\la \P(\X )\cup\{\emptyset\}, \subset\ra$. Since $C\neq \emptyset$ there is $A\in \L \setminus \{ \emptyset \}$
such that $A\subset C$ and, hence, $C\in \CP$. Thus $\L \cup \{ C\}$ is a chain in  $\CP\cup\{\emptyset\}$ bigger than
$\L$. A contradiction.

(b) follows from (a) and Fact \ref{T6024}.
\hfill $\Box$
\section{Maximal chains of copies of $\H_n$} \label{S4}
The graphs $\H _n$, $n\geq 3$, were constructed by Henson in \cite{henson}. By  \cite{henson}, $\H _n$ is  the unique (up to isomorphism)
countable ultrahomogeneous universal $\K _n$-free graph.

In order to cite a characterization of $\H _n$ which is more convenient for our construction,
we introduce the following notation.
If $\BG=\la G,\r \ra$ is a graph and $n\geq 3$ let $C_n(\BG )$ denote the set of all pairs
$\la H,K \ra$ of finite subsets of $G$ such that:

(C1) $K\subset H$ and

(C2) $K$ does not contain a copy of $\K _{n-1}$.

\noindent For $\la H,K \ra \in C_n(\BG )$, let $G^H_K$ denote
the set of all $v\in G\setminus H$ satisfying:

(S1) $ \{ v ,k \} \in \r$, for all $k\in K$ and

(S2) $ \{v, h \}  \notin \r$, for all $h\in H\setminus K$.

\noindent
The graphs $\H _n$ can be characterized in the following way.
\begin{fac}[Henson] \rm \label{T6036}
A countable graph $\BG=\la G,\r \ra$ is isomorphic to $\H _n$ iff it is $\K_n$-free and $G^H_K\neq \emptyset$, for each
$\la H,K \ra\in C_n(\BG )$.
\end{fac}
Now we prove (I) of Theorem \ref{T6021} for the graphs $\H _n$.
\begin{te}  \rm \label{T6001}
For each $n\geq3$ and each linear order $L$ the following conditions are equivalent:

(a) $L$ is isomorphic to a maximal chain in the poset $\langle \P(\H _n ) \cup \{ \emptyset \} , \subset \rangle $;

(b) $L$ is an ${\mathbb R}$-embeddable complete linear order with $0_L$ non-isolated;

(c) $L$ is isomorphic to a compact set $K \subset \mathbb R$ having the minimum non-isolated.
\end{te}
\dok
The equivalence (b) $\Leftrightarrow$ (c) is proved in Theorem 6 of \cite{K1} and (a) $\Rightarrow$ (b) follows from Theorem \ref{T5037}.

(b) $\Rightarrow$ (a).
We intend to use Theorem \ref{T6037}.
Let $\{ J_n' : n\in \o \}$ be a partition of the set $[0,1)\cap \Q$ into dense subsets of $[0,1)\cap \Q$. Let
$\mathbb Z$ denote the set of integers and let $J_n=\{ q+m : q \in J_n' \wedge m \in \mathbb Z \}$, for $n\in \o$.
Clearly, $\{ J_n : n\in \o \}$ is a partition of $\mathbb Q$ into dense subsets of $\mathbb Q$ and conditions (i) and (ii)
are satisfied.

Now we construct a copy of $\H _n$ with the domain $\Q$.
Let $\P$ be the set of $\K _n$-free graphs
$p=\la G_p , \r_p \ra$ such that $G_p \in [\Q ]^{<\o }$ and for each $a,b\in\Q$

(P1) $\{a,b\}\in \r _p \wedge \{a+1,b\}\in \r_p\Rightarrow b>a+1$,

(P2) $\{a,a-1\}\notin \r_p$.

\noindent
Let the relation $\leq$ on $\P$ be defined by
\begin{equation}\label{EQ6019}
p\leq q\Leftrightarrow G_p\supset G_q\ \wedge\ \r_p\cap [G_q]^2=\r_q.
\end{equation}
\begin{cla}\rm\label{T6008}
$\la \P , \leq \ra$ is a partial order.
\end{cla}
\dok
It is evident that the relation $\leq$ is reflexive and antisymmetric.
If $p\leq q\leq r$, then $G_r \subset G_q \subset G_p$ and $\r _r = \r _q \cap [G_r ]^2 = \r _p \cap [G_q ]^2\cap [G_r ]^2= \r _p \cap [G_r ]^2$.
\hfill $\Box$
\begin{cla}\rm\label{T6009}
The sets $\D _q =\{ p\in \P : q\in G_p\}$, $q\in \Q$, are dense in $\la \P , \leq \ra$.
\end{cla}
\dok
If $p=\la G_p , \r_p \ra \in \mathbb P\setminus \mathcal D_q$, then  $q\notin G_p$ and,
since $\{q,x\}\not\in \r _p$, for all $x\in G_p$,
$p_1=\la G_p\cup\{q\}, \r _p\ra$ is a $\K _n$-free graph and, clearly, satisfies (P1) and (P2). Thus $p_1\in \mathcal D_q$
and $p_1\leq p$.
\kdok
For $H\in [\Q ]^{<\o }$ let $m_H= \max H$.
\begin{cla}\rm\label{T6010}
For each $K\subset H\in[\Q]^{<\o}$ and each $m\in \mathbb N$, the set
$$\textstyle
\D^H_{K,m}=\Big\{p\in \P:H\subset G_p\ \wedge\ \Big(\la H,K \ra\in C_n(p) \Rightarrow \exists q\in J_0\cap (m_H,m_H+\frac{1}{m})
$$\\[-10mm]
$$\textstyle
\forall k\in K\; (\{q,k\}\in \r_p )  \land \forall h\in H\setminus K \; (\{q,h\}\not\in \r_p\ )\Big)\Big\}
$$
is dense in $\P$.
\end{cla}
\dok
Let $p_0\in \mathbb P$. By Claim \ref{T6009}
there is $p\in \mathbb P$ such that $p\leq p_0$ and $H\subset G_p$.

If $\la H,K \ra\notin C_n(p)$ then $p\in \D^H_{K,m}$ and we are done.

If $\la H,K \ra\in C_n(p)$, we take $q\in J_0\cap (m_H,m_H + \frac{1}{m})\setminus\bigcup_{a\in G_p}\{a,a-1,a+1\}$,
define
\begin{equation}\label{EQ6020}
p_1=\la G_p\cup\{q\}, \r _p \cup \{\{q,k\}:k\in K \}\ra .
\end{equation}
and first prove that $p_1\in \P$. Clearly $G_{p_1}\in [\Q]^{<\o}$ and we check that $p_1$ is $\K _n$-free.
Suppose that there is $F\in [G_{p_1}]^n$ such that $[F]^2\subset \r _{p_1}$. Since $p$ is $\K _n$-free
we have $q\in F$ and there are different $f_1,\dots,f_{n-1}\in G_p \cap F$
such that $\{q,f_i\}\in \r _{p_1}$, for $i\leq n-1$, which by (\ref{EQ6020}) implies
$\{f_1,\dots,f_{n-1}\}\subset K$.
Since $[F]^2\subset \r_{p_1}$, we have
$[\{f_1,\dots,f_{n-1}\}]^2\subset \r _{p}$.
But $\la H,K \ra\in C_n(p)$ implies that $K$ is $\K _{n-1}$-free.
A contradiction.

(P1) Suppose that for some $a,b\in \Q$
\begin{equation}\label{EQ6001}
\{a,b\}\in \r_{p_1}\ \wedge\ \{a+1,b\}\in \r_{p_1}\ \wedge b\leq a+1.
\end{equation}
Then, since $p\in\P$, at least one of the two pairs does
not belong to $\r _p$ and, hence, $q\in\{a,a+1,b\}$. So we have the following three cases.

$q=a$. Then by (\ref{EQ6001}) we have $b\neq q$ and, by (\ref{EQ6020}), $\{q+1,b\}\in \r_p$ which implies $q+1\in G_p$.
A contradiction to the choice of $q$.

$q=a+1$. Then by (\ref{EQ6001}) we have $b\neq q$ and, since $a\neq q$, by (\ref{EQ6020}) we have $\{a,b\}\in \r_p$
which implies $a\in G_p$. A contradiction to the choice of $q$.

$q=b$. Then by (\ref{EQ6001}) and (\ref{EQ6020}) we have $\{a,q\} , \{a+1,q\}\in \r_{p_1}\setminus\r_p$ which
implies $a,a+1\in K$. Since $q>m_H$ and $K\subset H$ we have
$q>a+1$ that is $b>a+1$. A contradiction again.

(P2) holds because $p\in \P$ and $q\not\in \bigcup_{a\in G_p}\{a,a-1,a+1\}$.

Thus $p_1\in \P$. Since $H\subset G_p \subset G_{p_1}$ and since, by (\ref{EQ6020}) we have
$\{q,k\}\in \r_{p_1} $, for all $k\in K$, and  $\{q,h\}\not\in \r_{p_1} $, for all $h\in H\setminus K$,
it follows that $p_1\in \D^H_{K,m}$.

Since $G_{p_1}\supset G_p$ and $\r _{p_1}\cap [G_p]^2=\r _p$, we have $p_1\leq p \leq p_0$.
\kdok
\noindent
By Fact \ref{T6023} there is a filter $\mathcal G$ in $\la \P , \leq \ra$ intersecting all sets $\D _q$, $q\in \Q$,
and $\D ^H_{K,m}$, for $K\subset H \in[\Q ]^{<\o}$ and $m\in \mathbb N$.
\begin{cla}\rm\label{T6028}
(a) $\bigcup _{p\in \mathcal G }G_p=\mathbb Q$;

(b) $\la \mathbb Q , \r \ra$ is a graph, where $ \r = \bigcup _{p\in \mathcal G}\r _p$, also $\{a,a-1\}\notin \r$, for all $a\in \Q$;

(c) $ \r \cap [G_p]^2=\r _p $, for each $p\in \mathcal G $;

(d) If $A\subset \Q$, $\r _A=\r \cap [A]^2$, $p\in \mathcal G$, and $H\subset A\cap G_p$, then $\r _A \cap [H]^2=\r _p\cap [H]^2$.
Thus if, in addition, $\la H,K \ra \in C_n (A, \r _A)$, then $\la H,K \ra \in C_n (p)$,

(e) $\la \mathbb Q , \r \ra$ is a $\K _n$-free graph.
\end{cla}
\dok
(a) For $q\in \mathbb Q$ let $p_0\in \mathcal G \cap \mathcal D_q$. Then $q\in G_{p_0} \subset \bigcup_{p\in \mathcal G}G_p$.

(b) By the definition of $\P$ we have $\{a,a-1\}\notin \r _p \subset [\Q ]^2$, for all $p\in \P$.

(c) The inclusion ``$\supset$" is evident. If $\{ a,b \} \in \r \cap [G_p]^2$,
then there is $p_1\in \mathcal G$ such that $\{ a,b \} \in \r _{p_1}$ and, since
$\mathcal G $ is a filter, there is $p_2\in \mathcal G$ such that
$p_2\leq p ,p_1 $. By the definition of $\leq$ we have $\r _{p_1}\subset \r _{p_2}$, which implies $\{ a,b \} \in \r _{p_2}$
and $\{ a,b \} \in \r _{p_2}\cap [G_p]^2=\r _p$.

(d) By (c) we have $\r _A \cap [H]^2 = \r \cap [A]^2 \cap [H]^2= \r \cap [H]^2= \r \cap [G_p]^2 \cap [H]^2= \r _p \cap [H]^2$.
If $\la H,K \ra \in C_n (A,\r _A)$, then $K$ is $\K _{n-1}$-free in $\la A, \r _A \ra$ and, since $\r _A \cap [K]^2=\r _p\cap [K]^2$, $K$ is $\K _{n-1}$-free in $p$ as well. Thus $\la H,K \ra \in C_n (p)$.

(e) Suppose that $\la A , \r_A \ra$ is a copy of $\K _n$ and let $p_q\in \mathcal G \cap \D _q$, $q\in A$. Since $\mathcal G$ is a filter
there is $p\in \mathcal G$ such that $p\leq p_q$, for all $q\in A$, and, hence, $A\subset G_p$, which by (d)
implies $\r _A =\r _p \cap [A]^2$. But this is impossible since $p$ is $\K _n$-free.
\kdok
Now we show that conditions (iii) and (iv) of Theorem \ref{T6037} are satisfied.

(iii) Let $x\in \R \cup \{ \infty \}$ and $J_0 \cap (-\infty ,x) \subset  A\subset \Q \cap (-\infty , x) $. We show that $\la A ,\r _A \ra \cong \H _n$.
By Claim \ref{T6028}(e) $\langle A,\r _A \rangle$ is $\K _n$-free.
Let $\la H,K \ra \in C_n (A, \r _A)$. Since $m_H \in H\subset A$ we have $m_H<x$
and there is $m\in \N$ satisfying $m_H+\frac{1}{m}<x$.
Let $p\in \mathcal G \cap \D^H_{K,m}$. Then $K\subset H\subset G_p$ and,
by Claim \ref{T6028}(d), $\la H,K\ra\in C_n (p)$.
Thus there is $q\in J_0 \cap (m_H,m_H+\frac{1}{m})\subset J_0\cap (-\infty,x)\subset A$ such that
$\{ q,k \}\in \r _p \subset \r$, which implies $\{ q,k \}\in \r _A$, for all $k\in K$,
and that $\{ q,h \}\not\in \r _p$, which implies $\{ q,h \}\not\in \r $, for all $h\in H$.
Thus $q\in A^H_K$. By Fact \ref{T6036} we have $\la A ,\r _A \ra \cong \H _n $.

(iv) Let $q\in J_0$ and $J_0 \cap (-\infty ,q] \subset  C\subset \Q \cap (-\infty , q]$. We prove that $ \la C ,\r _C \ra \not\cong \H _n$.
Since $q\in J_0$ by the construction of $J_0$ we have $q-1 \in J_0$ and, by the assumption,
$H=\{ q-1, q \}\subset C$. By Claim \ref{T6028}(b) we have $\{ q-1, q \}\not\in \r$ which implies that
$H$ is $\K _{n-1}$-free and, hence, $\la H,H \ra \in C_n (C, \r _C)$.
Suppose that $b\in C^H_H$. Then $\{ q-1, b \}, \{q,b \} \in \r$ and, since $\mathcal G$ is a filter,
$\{ q-1, b \}, \{q,b \} \in \r _p$, for some $p\in \mathcal G$. By (P1) we have $b>q$, which is impossible since
$q=\max C$. Thus $C^H_H=\emptyset$ and by Fact \ref{T6036} we have $\la C ,\r _C\ra \not\cong \H _n $.
\begin{cla}\rm\label{T6038}
The family
$\CP=\Big\{ \Q\setminus\bigcup_{n\in\Z}F_n :\forall n\in \Z \;\; F_n \in \Big[[n,n+1)\cap\Q\Big]^{<\omega}\Big\}$ is
a positive family in $P(\Q)$ satisfying $\bigcap \CP=\emptyset$ and $\CP\subset \P(\Q,\r )$.
\end{cla}
\dok
It is easy to check (P1)-(P4). Since $\Q \setminus \{ q\}\in \CP$, for each $q\in \Q$, we have $\bigcap \CP=\emptyset$.
Let $A=\Q\setminus\bigcup_{n\in\Z}F_n\in \CP$, $\la H,K \ra\in C_n(A, \r _A)$ and $m_H=\max H\in [n_0,n_0 +1)\cap \Q$.
Since $|F_{n_0}|<\o$ and $m_H \in A \subset \Q\setminus F_{n_0}$
there is $m\in \N$ such that $(m_H,m_H+\frac{1}{m})\cap \Q\subset A$.
Let $p\in \mathcal G\mathcal D^H_{K,m}$. Then $H\subset G_p$ and, by Claim \ref{T6028}(d), $\la H,K \ra\in C_n(p)$.
Hence there is $q\in J_0 \cap (m_H,m_H+\frac{1}{m})\subset A$
such that

- for each $k\in K $ we have $\{ q , k \}\in \r _p$ which, since $\{ q , k \}\subset A\cap G_p$, by Claim \ref{T6028}(d) implies
$\{ q , k \}\in \r _A$;

- for each $h\in H\setminus K$ we have $\{ q , h \}\not\in \r _p$, which by Claim \ref{T6028}(c) implies $\{ q , h \}\not\in \r $ and, hence,
$\{ q , h \}\not\in \r _A$.

\noindent
Thus $q\in A^H_K$. By Fact \ref{T6036} we have $\la A, \r _A \ra \cong \H _n \cong \la \Q ,\r \ra$ and, hence, $A\in \P(\Q,\r )$.
\kdok
Now (b) $\Rightarrow $ (a) of Theorem \ref{T6001} for countable $L$ follows from Claim \ref{T6038} and Theorem \ref{T6025}(b).
Thus condition (v) of Theorem \ref{T6037} is satisfied and, by (B) of Theorem \ref{T6037}, (b) $\Rightarrow $ (a) of Theorem \ref{T6001} is true for uncountable $L$.
\hfill $\Box$
\section{Maximal chains of copies of $\BG _{\mu \nu}$}\label{S5}
\begin{te}\rm\label{T6040}
If $\mu$ and $\nu$ are cardinals satisfying $\mu \nu =\o$, then for each linear order $L$ the
following conditions are equivalent:

(a) $L$ is isomorphic to a maximal chain in the poset $\langle \P(\BG_{\mu \nu}) \cup \{ \emptyset \} , \subset \rangle $;

(b) $L$ is an ${\mathbb R}$-embeddable Boolean linear order with $0_L$ non-isolated;

(c) $L$ is isomorphic to a compact nowhere dense set $K \subset \mathbb R$ having the minimum non-isolated.
\end{te}
\dok
Clearly, concerning the values of $\mu$ and $\nu$ we have three cases.

I. $\BG _{\o n}=\bigcup _{i\in \o } \BG _i$, where $n\in \N$ and $\BG _i =\la G_i , [G_i]^2\ra$, $i\in \o$, are disjoint copies of $\K _n$.
Then, clearly $\P (\BG _{\o n}) = \{ \bigcup_{i\in A}G_i : A\in [\omega]^{\omega}\}$ and, hence,
$\la \P (\BG _{\o n}) \cup \{ \emptyset \}, \subset \ra \cong \la [\o]^\o \cup \{ \emptyset \}, \subset\ra$.
Since $[\o ]^\o$ is a positive family in $P(\o )$ the statement follows from Fact \ref{T6022}.

II. $\BG _{m \o }=\bigcup _{i<m } \BG _i$, where $m\in \N$ and $\BG _i =\la G_i , [G_i]^2\ra$, $i<m$, are disjoint copies of $\K _\o$.
Then, since each copy of $\BG _{m \o }$ must have $m$ components of size $\o$, we have
$\P (\BG _{m \o }) =\{  \bigcup _{i<m}A_i : \forall i<m \; A_i \in [G_i]^\o \}$ and it is easy to see that
$\P (\BG _{m \o })$ is a positive family in $P(G)$ so we apply Fact \ref{T6022} again.

III. $\BG _{\o \o }=\bigcup _{i<\o } \BG _i$, where $\BG _i =\la G_i , [G_i]^2\ra$, $i<\o$, are disjoint copies of $\K _\o$.
The equivalence (b) $\Leftrightarrow$ (c) is a part of Fact \ref{T6022}

(a) $\Rightarrow$ (b). If $\L$ is a maximal chain in $\langle \P(\BG _{\o \o } ) \cup \{ \emptyset \} , \subset \rangle $, then,
by Theorem \ref{T5037}, it is an ${\mathbb R}$-embeddable complete linear order with $0_\L$ non-isolated and we prove that it has
dense jumps. Let $G=\bigcup _{i<\o}G_i$. Since each copy of $\BG _{\o \o }$ must have $\o$ components of size $\o$, we have
\begin{equation}\textstyle \label{EQ6024}
\P (\BG _{\o \o })=\{ \bigcup _{i\in S}A_i : S\in [\o ]^\o \land \forall i\in S \; A_i \in [G_i ]^\o \}
\end{equation}
and, for $A= \bigcup _{i\in S}A_i \in \P (\BG _{\o \o })$ we will write $S=\supp A$.

Let $A,B\in \L \setminus \{\emptyset \}$, where $A\varsubsetneq B$.
\begin{cla}\rm \label{T6035}
There is $C\in \L$ satisfying $A\subset C\subset B$ and such that $C\cap G_i \varsubsetneq B\cap G_i$, for some $i\in \supp C$.
\end{cla}
\dok
Suppose that for each $C\in \L \cap [A,B]$ we have: $C\cap G_i = B\cap G_i$, for all $i\in \supp C$. Then, since $A\varsubsetneq B$, we have
$\supp A\varsubsetneq \supp B$ and we choose $i\in \supp B\setminus \supp A$.
Clearly, for the sets $\L ^- =\{ C\in \L: i\not\in \supp C \}$ and $\L ^+ =\{ C\in \L: i\in \supp C \}$
we have $\L =\L ^- \cup \L^+$ and $C_1\varsubsetneq C_2$, for each $C_1 \in \L ^-$ and $C_2\in \L ^+$.
By Theorem \ref{T5037}  we have $C^-=\bigcup \L ^- \in \P (\BG _{\o \o })$ and, since $\L^- \vartriangleleft \L ^+$, by the maximality of $\L$ we have $C^-\in \L$.
Clearly $i\not\in \supp C^- $, which implies $C^- =\max \L ^-$.
Let $C^+=C^- \cup (B\cap G_i)$. By (\ref{EQ6024}) we have $C^+\in \P (\BG _{\o \o })$.
For $C\in \L ^+$ we have $i\in \supp C$ and, by the assumption,
$C\cap G_i = B\cap G_i$, which implies $C^+\subset C$. Thus, by the maximality of $\L$, $C^+ \in \L$, and, moreover, $C^+=\min \L ^+$.
Let $a\in B\cap G_i$. Then $C=C^- \cup (B\cap G_i \setminus \{ a\})\in \P (\BG _{\o \o })$ and $C^- \varsubsetneq C \varsubsetneq C^+$, which implies that
$\L$ is not a maximal chain in $\P (\BG _{\o \o })$. A contradiction.
\kdok
Let $C_0\in \L$ and $i_0\in \supp C_0$ be the objects provided by Claim \ref{T6035}. Let $a\in (B\setminus C_0)\cap G_{i_0}$,
$\L ^- =\{ C\in \L: a\not\in C \}$ and $\L ^+ =\{ C\in \L: a\in C \}$.
Then we have $\L =\L ^- \cup \L^+$, $C_0\in \L ^-$ and $C_1\varsubsetneq C_2$, for each $C_1 \in \L ^-$ and $C_2\in \L ^+$.
By Theorem \ref{T5037} we have $C^-=\bigcup \L ^- \in \P (\BG _{\o \o })$ and, by the maximality of $\L$, $C^- \in \L$.
Since $a\not\in C^-$ we have $C^- =\max \L ^-$, which implies $C_0 \subset C^-$ and, hence, $i_0\in \supp C^-$.
Thus, by (\ref{EQ6024}),
$C^+ =C^- \cup \{ a \}\in \P (\BG _{\o \o })$. For $C\in \L ^+$ we have $C^+ \subset C$ and, by the maximality of $\L$, $C^+ \in \L$, in fact $C^+ =\min \L^+$.
Clearly the pair $\la C^- ,C^+\ra$ is a jump in $\L$. Since $A\subset C_0$ and $B\in \L ^+$ we have
$A\subset C^- \subset C^+ \subset B$.
Thus, $\L \setminus \{ \emptyset \}$ has dense jumps and, since $0_\L$ is non-isolated, the same holds for $\L$.

(b) $\Rightarrow$ (a). Clearly, $\CP =\{ \bigcup _{i\in \o}A_i : \forall i\in \o \; A_i \in [G_i ]^\o \}$ is a positive family contained in
$\P (\BG _{\o \o })$ and $\bigcap \CP =\emptyset$. Now the statement follows from Theorem \ref{T6025}(a).
\hfill $\Box$

\end{document}